\documentclass{amsart}
\usepackage{amssymb}
\title{On monotonicity of F-blowup sequences}
\author{Takehiko Yasuda}
\address{Department of Mathematics and Computer Science, 
Faculty of Science, Kagoshima University, 1-21-35 Korimoto, Kagoshima 890-0065, Japan}
\email{yasuda@sci.kagoshima-u.ac.jp}
\thanks{I thank
 Shunsuke Takagi, Kei-ichi Watanabe and Ken-ichi Yoshida
 for helpful comments concerning F-singularities.
 I also thank the referee for his/her helpful comments and raising a few interesting 
 questions, which are answered in 
 Propositions \ref{prop-not-preserve} and \ref{prop-F-pure-weakly-normal}.}

\theoremstyle{plain}
\newtheorem{thm}{Theorem}[section]

\newtheorem{prop}[thm]{Proposition}
\newtheorem{cor}[thm]{Corollary}
\newtheorem{lem}[thm]{Lemma}

 \theoremstyle{definition}
\newtheorem{defn}[thm]{Definition}
\newtheorem{ex}[thm]{Example}

\theoremstyle{remark}

\def\AA{\mathbb A}

\def\RR{\mathbb R}

\def\ZZ{\mathbb Z}

\def\Zp{\mathbb{Z}_{>0}}
\def\Znn{\mathbb{Z}_{\ge 0}}

\def\cF{\mathcal{F}}
\def\cG{\mathcal{G}}

\def\cM{\mathcal{M}}

\def\cO{\mathcal{O}}

\def\fz{\mathfrak z}

\def\11{\mathbf{1}}

\def\id{\mathrm{id}}

\DeclareMathOperator{\Spec}{Spec}

\DeclareMathOperator{\Hom}{Hom}

\DeclareMathOperator{\FB}{FB}

\DeclareMathOperator{\Hilb}{Hilb}
\DeclareMathOperator{\Quot}{Quot}

\def\GHilb{\mathrm{Hilb}^{G}}

\begin{document}

\maketitle

\begin{abstract}
For each variety in positive characteristic,
there is a series of canonically defined blowups, called F-blowups. 
We  are interested in the question of whether the $e+1$-th
blowup dominates the $e$-th, locally or globally.
It is shown that the answer is affirmative (globally for any $e$) 
when the given variety is F-pure. 
As a corollary, we obtain some result on the stability
of the sequence of F-blowups. 
We also give a sufficient condition for local domination.  
\end{abstract}

\section{Introduction}

The F-blowup introduced in \cite{Yasuda:math0706.2700} 
is an interesting notion which relates for instance to
 the desingularization problem, the $G$-Hilbert
scheme and Gr\"obner bases. 
The study of it has just started and 
there remain various problems.
Among them, it seems important to understand the
behavior of the sequence consisting of F-blowups. 

Consider a variety $X$ in positive characteristic, that is,
a separated integral scheme of finite type over an algebraically
closed field $k$ of characteristic $p > 0$. 
Let 
\begin{equation*}
F_{e} : X_{e} \to X , \ e =0,1,2, \dots
\end{equation*}
be the $e$-times iteration of the $k$-linear Frobenius. 
Then for each smooth (closed)  point $x \in X$, the fiber $F_{e}^{-1}(x)$
is a fat point of $X_{e}$ of length $p^{e \cdot \dim X}$.
It is considered as a reduced point of the Hilbert scheme 
of $0$-dimensional subschemes of $X_{e}$ of this length:
 $F_{e}^{-1}(x) \in \Hilb _{p^{e \cdot \dim X}} (X_{e})$.

\begin{defn}
We define the \emph{$e$-th F-blowup}
of $X$, $\FB_{e} (X)$, as the closure of 
the subset
\[
 \{ F_{e}^{-1}(x) |x \in X(k) \ \text{smooth} \} \subset \Hilb_{p^{e \cdot \dim X}} (X_{e}).
\]
\end{defn}

This is indeed a blowup of $X$,
that is, birational and projective over $X$ (Proposition \ref{prop-relative}).

It is natural to ask if  
$\FB_{e+1} X $ dominates $ \FB_{e} X$,
that is, if the natural birational map 
\[
\rho _{e} :\FB_{e+1} X\dasharrow  \FB_{e}X 
\]
has no indeterminacy. 
When this holds for all $e$, we shall say the F-blowup sequence
is \emph{monotone}.   
The answer is generally negative (Example \ref{ex}).
One of our main theorems provides a sufficient condition for the monotonicity:

\begin{thm}\label{thm-F-pure}
Suppose that $X$ is F-pure,
that is,  the natural morphism $\cO_{X} \to (F_{1})_{*} \cO_{X}$ locally splits
as a morphism of $\cO_{X}$-modules. 
Then the F-blowup sequence of $X$ 
is monotone.
\end{thm}

The notion of F-purity was introduced by Hochster and Roberts \cite{MR0417172}
and is now one of important classes of F-singularities (see  \cite{smith-guide}
and the references given there).

We can consider also the local version of the above question:
``Is $\rho_{e}$ 
defined at a given point of $\FB_{e +1} X$?'' 
We give a sufficient condition for this too;

\begin{thm}\label{thm-local-condition}
Let $Z \in \FB_{e +1} X$ be a closed point, which is identified with a fat point of $X_{e+1}$.
Suppose  that 
\begin{quotation}
(*) the scheme-theoretic image $\bar Z$ of $Z$
by the natural morphism $X_{e+1} \to X_{e}$ belongs to $\FB_{e} X$.
\end{quotation}
Then $\rho_{e}$ is
defined at $Z$ and $\rho_{e}(Z)= \bar Z$.
\end{thm}

Condition (*) means that
$Z \in \FB_{e +1} X$ has a natural candidate $\bar Z$ of the image
in $\FB_{e} X$. 
In fact, Theorem \ref{thm-local-condition} is a generalization of 
Theorem \ref{thm-F-pure}, since (*)
always holds if $X$ is F-pure (Proposition \ref{prop-onto}).

We saw in \cite{Yasuda:math0706.2700} that in some cases,
the F-blowup sequence is bounded, that is, all
F-blowups of $X$ are dominated by a single blowup of $X$.
(At this point I do not know of any example where the sequence is unbounded.
See Example \ref{examples}.)
When both the boundedness and monotonicity hold,
the sequence stabilizes (Lemma \ref{lem-stability}).

It is also natural to ask what properties of variety 
are preserved by F-blowups. 
We obtain the following result on this issue:

\begin{prop}\label{prop-not-preserve}
There exist an F-pure (resp.\ normal, weakly normal) variety $X$ and $e \in \Zp$
such that $\FB_{e} X$ is not F-pure  (resp.\ normal, weakly normal).
\end{prop}

We use the  toric  geometry 
in order to construct examples and prove the last proposition.
For this purpose, we show that a toric variety is
F-pure if and only if it is weakly normal (Proposition \ref{prop-F-pure-weakly-normal}). 
A similar result has been obtained by Bruns, Li and R\"omer \cite[Proposition 6.2]{MR2236607}.

\subsection*{Outline of the paper}
In Section \ref{sec-prelim}, we recall some basic facts on F-blowup
from \cite{Yasuda:math0706.2700}.
Section \ref{sec-proof} is devoted to the proof of Theorem \ref{thm-F-pure}.
In Section \ref{sec-local-condition}, we prove that 
Condition (*) holds whenever $X $ is F-pure,
and Theorem \ref{thm-local-condition}.
In Section \ref{sec-toric}, we use the toric geometry
to give some examples of F-blowups.
In Section \ref{sec-stability}, we discuss when the F-blowup sequence stabilizes.
In Section \ref{sec-not-preserve}, we prove Proposition \ref{prop-not-preserve}
by using the toric geometry and the nonnormal $G$-Hilbert scheme 
found by Craw, Maclagan and Thomas \cite{MR2356842}.

\subsection*{Conventions}
Throughout the paper, we work over an algebraically closed field $k$
of characteristic $p >0$.
A variety means a separated integral scheme of finite type
over $k$.
A point of a variety always means a closed point.

\section{Preliminaries}\label{sec-prelim}

In this section, we set up notation and recall some
results from \cite{Yasuda:math0706.2700}.\footnote{
In \cite{Yasuda:math0706.2700}, the Frobenius map
$R^{q} \hookrightarrow R$,
 rather than $R \hookrightarrow R^{1/q}$,
is considered. 
This causes slight notational differences.}

We continue to write $X$ for a given variety over $k$.
All our problems are local on $X$.
So we may suppose $X$ is affine; $X = \Spec R$. 
Let $e \in \Znn$ and $q :=p^{e}$.
Then we may identify $X_{e} = \Spec R^{1/q} $
and  then $F_{e} :X_{e} \to X$ corresponds to the inclusion
map $R \hookrightarrow R^{1/q}$. 
We also have $(F_{e})_{*} \cO_{X} = \cO_{X}^{1/q} $ 
and $X_{e} = \mathcal{S} pec_{X} \, \cO_{X}^{1/q}$.

The F-blowup can be constructed also with the \emph{relative}
Hilbert scheme or the Quot scheme:

\begin{prop}\label{prop-relative}\cite[Proposition 2.4]{Yasuda:math0706.2700}
The F-blowup $\FB_{e} (X)  $ is canonically isomorphic to the irreducible component of 
 $\Hilb_{q^{d}} (X_{e}/X) $
that dominates $X$, and also to that of $ \Quot _{q^{d}} (\cO_{X}^{1/q})$. 
\end{prop}

Moreover the proof of \cite[Proposition 2.4]{Yasuda:math0706.2700}
shows that the isomorphism is the restriction of the natural 
morphism $\Hilb_{q^{d}}(X_{e}/X) \to \Hilb_{q^{d}}(X_{e})$.
It follows that  each point $Z \in \FB_{e} (X)$ is included in
the fiber $F_{e}^{-1}(x)$ for some reduced point $x \in X$.
Namely the scheme-theoretic image $F_{e} (Z) \subset X$ is a reduced point.
Then the $X$-scheme structure of $\FB_{e} (X)$
is given by the map
\begin{equation*}\label{eq-pi_{e}}
\pi_{e}: \FB_{e} (X) \to X,\ Z \mapsto \pi_{e}(Z):= F_{e}(Z).
\end{equation*}
This is projective and is an isomorphism exactly
over the smooth locus of $X$ \cite[Corollary 2.5]{Yasuda:math0706.2700}. 

Being an irreducible component of the Quot scheme, 
 $\FB_{e} (X) $ has the following universal property:
For a blowup $f:Y \to X$ and a coherent $\cO_{X}$-module $\cF$,
define the \emph{torsion-free pull-back}
$f^{\bigstar} \cF$ the quotient of the usual pull-back $f^{*} \cF$
by the subsheaf of torsions. 
Then $\pi_{e}^{\bigstar} \cO_{X}^{1/q} $ is flat  or equivalently locally free.
Moreover if for a blowup $f:Y \to X$, $f^{\bigstar} \cO_{X}^{1/q} $ is flat,
then $f$ factors through $\FB_{e} X$.

More generally, if $\cG$ is a coherent sheaf on $X$
and if it is generically locally free of rank $r$, then 
its universal (birational) flattening is constructed as
the irreducible component of $\Quot_{r}(\cG)$ dominating $X$.
See for instance \cite{MR1218672,1087.14011} 
for  studies on the universal flattening of a general coherent module.

\section{Proof of Theorem \ref{thm-F-pure}}\label{sec-proof}

We may suppose that $X$ is affine.
Then for each $e$, we have an isomorphism of $\cO_{X}$-modules
\[
  \cO^{1/p^{e+1}}_{X} \cong \cO^{1/p^{e}}_{X} \oplus \cM _{e}
\]
for some $\cO_{X}$-module $\cM_{e}$.
Then the torsion-free pull-back by  $\pi_{e+1}$ 
\[
\pi_{e+1}^{\bigstar} \cO^{1/p^{e+1}}_{X}\cong \pi_{e+1}^{\bigstar}\cO^{1/p^{e}}_{X}\oplus \pi_{e+1}^{\bigstar} \cM_{e}
\]
is flat and locally free. 
From the characterization of flat module as a  summand of a free module \cite[Corollary 6.6]{MR1322960},
 $ \pi_{e+1}^{\bigstar}\cO^{1/p^{e}}_{X}$  is flat.
 From the universality of $\FB_{e} X$, we have a natural morphism
$ \FB_{e +1} X \to \FB_{e} X$. We have proved the theorem.

\section{On local domination by
$\FB_{e+1} X $ over $\FB_{e} X$}\label{sec-local-condition}

\begin{prop}\label{prop-onto}
Suppose that $X$ is F-pure. Then Condition 
(*) in Theorem \ref{thm-local-condition} holds
for every $e \ge 0$ and every $Z \in  \FB_{e+1}X$.
\end{prop}

\begin{proof}
The identity map of $\pi_{e+1}^{\bigstar}\cO^{1/p^{e}}_{X}$
can be
factored as
\[
\pi_{e+1}^{\bigstar}\cO^{1/p^{e}}_{X} \to \pi_{e+1}^{\bigstar}\cO^{1/p^{e+1}}_{X} \to
 \pi_{e+1}^{\bigstar}\cO^{1/p^{e+1}}_{X} /\pi_{e+1}^{\bigstar} \cM_{e} \cong \pi_{e+1}^{\bigstar}\cO^{1/p^{e}}_{X} .
\]
Taking the fibers of these sheaves at $Z$,
we obtain 
\[
\id_{k[Z']} : k[Z'] \to k[Z] \to k[Z'] .
\]
Here $Z' \in \FB_{e} X$ is the image of $Z$ by
the natural map $\FB_{e+1} X \to \FB_{e} X$,
which exists from Theorem \ref{thm-F-pure},
and $k[Z]$ and $k[Z']$
are the coordinate rings of fat points $Z \subset X_{e+1}$
and $Z' \subset X_{e}$ respectively. 
Hence the map $k[Z'] \to k[Z]$, which is the ring homomorphism
defining the natural morphism $Z \to Z'$,
is injective. This shows that $Z' = \bar Z$ 
and the proposition follows.
\end{proof}

\begin{proof}[Proof of Theorem \ref{thm-local-condition}]
We write  $ Z_{e+1} := Z $ and $ Z_{e} := \bar Z$.
Let $G \subset \FB_{e+1} (X) \times_{k} \FB_{e}(X)$ be the closure of the graph
of $\rho_{e} : \FB_{e+1} (X) \dasharrow \FB_{e}(X)$ and $\psi_{i}:G \to \FB_{i}(X)$, $i=e,e+1$,
the projections.
We have to show that  $\psi_{e+1}$
is an isomorphism around $a:= (Z_{e+1},Z_{e}) \in G$.

We shall first show that set-theoretically $\psi_{e+1}^{-1}(Z_{e+1})=\{a\}$.
For $i=e,e+1$, let $W_{i} \subset G \times _{k} X_{i}$ be
the family of fat points over $G$.
More precisely, this is
the pull-back of the universal family over $\FB_{i}(X)$,
which is a closed subscheme of $\FB_{i}(X) \times_{k} X_{i}$,
by the projection $\psi_{i}$.
Then $W_{i}$ is isomorphic to the associated reduced scheme of
$G \times_{X} X_{i}$. 
In other words, $\cO_{W_{i}}$ is identified with the torsion-free pull-back
of $\cO_{X_{i}}=\cO_{X}^{1/p^{i}}$ by the natural map $G \to X$.
In particular, $W_{i}$ is reduced. 
Hence  $W_{e}$ is the scheme-theoretic image of $W_{e+1}$ 
by the natural morphism $G \times_{k} X_{e+1} \to G \times_{k} X_{e}$.
If $b=(Y_{e+1},Y_{e}) \in G$, then   the fiber of $W_{e+1} \to G$ (resp.\ 
$W_{e} \to G$) over $b$ is $Y_{e+1}$ (resp.\ $Y_{e}$).
It follows that the scheme-theoretic image $\bar Y_{e+1}$ of $Y_{e+1}$ in $X_{e}$
is included in $Y_{e}$.
Now if $Y_{e+1} =Z_{e+1}$, then $Z_{e}:= \bar Y_{e+1}  \subset Y_{e}$.
But by assumption both $Z_{e}$ and $Y_{e}$ have length $p^{ed}$. Hence $Z_{e} =Y_{e}$.
This shows that $\psi_{e+1}^{-1}(Z_{e+1})=\{a\}$.

Let $R$ be the coordinate ring of $X$ as before,
 $\fz _{i} \subset R^{1/p^{i}}$
the defining ideals of $Z_{i}$ and
\[
\phi_{i}: T _{a} G \to T_{Z_{i}} \FB_{i}(X) \hookrightarrow \Hom( \fz_{i} , R^{1/p^{i}}/\fz_{i}) 
\]
the maps of Zariski tangent spaces (for the identification of the tangent space to the Hilbert scheme with $\Hom( \fz_{i} , R^{1/p^{i}}/\fz_{i})$, see for instance \cite[Proof of Theorem VI-29]{MR1730819}).
To show that $\psi_{e+1}$ is an isomorphism around $a$,
 it is enough to show that $\phi_{e+1}$ 
is injective. 
Take   $0 \ne v \in T_{a} G$. If $\phi_{e}(v)=0$, then $\phi_{e+1}(v) \ne 0$.
So we may suppose that $\phi_{e}(v)\ne 0$.
Let 
\[
W_{e}^{v} \subset \Spec R^{1/p^{e}}[t]/ (t^{2}) \text{ and }
W_{e+1}^{v} \subset \Spec R^{1/p^{e+1}}[t]/(t^{2})
\]
 be the pull back of $W_{e}$ and $W_{e+1}$
by our tangent vector
\[
v:\Spec k[t]/(t^{2}) \to G.
\]
Since $\phi_{e} (v) \ne 0$, the defining ideal of $W_{e}^{v}$
does contain an element of the form $f + g t$, $f \in \fz_{e}$, $g \in R^{1/p^{e}} \setminus \fz_{e}$
so that $\phi_{e}(v) \in \Hom (\fz_{e}, R^{1/p^{e}} / \fz_{e})$
maps $f$ to the class of $g$ modulo $\fz_{e}$, which is nonzero.

Such an element $f+gt$ is also contained in the defining ideal of $W_{e+1}^{v}$
 and $\phi_{e+1}(v)$ maps $f$ to $g$ modulo $\fz_{e+1}$.
Since by assumption $\fz_{e} = \fz_{e+1} \cap R^{1/p^{e}}$, we have 
$g \notin \fz_{e+1}$. Hence $\phi_{e+1}(v) \ne 0$
and $\phi_{e+1}$ is injective, which completes the proof. 
\end{proof}

\section{The toric case}\label{sec-toric}

Let $M = \ZZ^{d}$ be a free abelian group of rank $d$
and $A \subset M$ a finitely generated submonoid
which generates $M$ as a group. 
We associate to $A$ and $M$ the monoid algebras
 $k[A] \subset k[M] =\bigoplus_{m \in M} k \cdot x^{m}$
and the affine toric varieties $X := \Spec k[A] \supset T := \Spec k[M]$.
We shall make an additional assumption that
  $A$ contains no nontrivial group or equivalently
the cone $A_{\RR} \subset M_{\RR}$ spanned by $A$
has a vertex.
This involves no loss of generality.\footnote{
Conversely suppose that  $A$ contains a nontrivial group.
Let $B \subset A$ be the maximal group
and let  $a_{1} , \dots, a_{m},b_{1}, \dots,b_{n } $ be 
generators of $A$ such that $a_{i} \notin B$ and $b_{i} \in B$.
Then there exists a subset of $\{b_{1}, \dots, b_{n}\}$,
say $\{b_{1}, \dots, b_{l}\}$, $l \le n$,
which generates a monoid containing no nontrivial group
but still generates $B$ as a group.
Let  $A'$ be the monoid generated by $a_{1} , \dots, a_{m},b_{1}, \dots,b_{l } $,
which contains no nontrivial group. 
Then $k[A]$ is a localization of $k[A']$ by an element. 
Indeed if we put $b:= \sum_{i=1}^{l} b_{i}$,
then  $k[A ] = k[A']_{x^{b}}$. Thus the toric variety associated to $A$
is an open subvariety of the one associated to $A'$.}

Let $A_{\RR}^{\vee} \subset M _{\RR}^{\vee}$ be the dual cone 
of $A_{\RR}$, which is $d$-dimensional since $A_{\RR}$ is strongly
convex. The  $F$-blowup $\FB_{e} X$ is a (possibly nonnormal) 
toric variety and determines
a fan $\Delta _{e}$ which is a subdivision of $A_{\RR}^{\vee}$.
For each $d$-dimensional cone $\sigma \in \Delta_{e}$,
there exists a corresponding affine toric open subvariety  $U_{\sigma} 
\subset X$. 
By the inclusion $T \subset U_{\sigma}$, 
the coordinate rings of each $U_{\sigma}$ is naturally embedded in $k[M]$.
It is expressed as follows:
Fix a $d$-dimensional $\sigma \in \Delta_{e}$
and an interior point $w \in \sigma$.
Then put 
\begin{equation*}\label{def-B}
 B_{\sigma}  := \{ m \in \frac{1}{q} A  |  \exists m' \in \frac{1}{q} A , \,   m-m' \in M, \, \langle m,  w \rangle > \langle m', w \rangle \}, 
\end{equation*}
and
\begin{multline*} \label{def-C}
C_{\sigma} := \{ m-m' | 
  m \in \frac{1}{q} A,\, m' \in \frac{1}{q} A \setminus B_{\sigma},\,  
 m-m' \in M , \, \langle m,  w \rangle > \langle m', w \rangle  \}.
\end{multline*}

\begin{thm}\cite[Proposition 3.8]{Yasuda:math0706.2700}
The coordinate ring of $U_{\sigma}$ is generated by 
$x^{c}$, $c \in C_{\sigma}$ as a $k$-algebra.
\end{thm}

Now we recall the notion of weak normality 
in the sense of Andreotti and Bombieri \cite{MR0266923}:

\begin{defn}
An affine variety $\Spec R$ with function field $K$ 
is said to be \emph{weakly normal} if $R=R^{1/p} \cap K$.
\end{defn}

We easily see that the monoid algebra $k[A]$ is weakly normal if and only if 
$A= \frac{1}{p} A  \cap M$.

It has been known to experts
 that the F-purity implies the weak normality (for example, see \cite[Proposition 1.2.5]{MR2107324}).
For the monoid algebra, the converse is also true:

\begin{prop}\label{prop-F-pure-weakly-normal}
The ring $k[A]$ is F-pure if and only if it is weakly normal.
\end{prop}

\begin{proof}
Although I do not know of any reference, maybe this result is known to experts.
Bruns, Li and R\"omer \cite[Proposition 6.2]{MR2236607}
have proved a similar result. 
Suppose that $k[A]$ is weakly normal, so $A= \frac{1}{p}A  \cap M$.
We define a $k$-linear map $\phi: k[\frac{1}{p}A] \to k[A]$ by
\[
\phi (x ^{m}) = 
\begin{cases}
x^{m} &  (m \in A) \\
0 & (m \notin A) . 
\end{cases}
\]
We claim that it is a $k[A]$-module homomorphism
and hence the inclusion map $k[A] \hookrightarrow k[\frac{1}{p}A]$ splits. 
To see this,  it is enough to show that for any $m \in \frac{1}{p} A$ and $n \in A$,
\begin{equation}\label{eq-phi}
\phi(x^{m +n}) = x^{n}\phi(x^{m}) .
\end{equation}
When $m \in A$, this is obvious. 
If $m \notin A$, then $m + n \notin A$.
(Conversely, if $m + n \in A$, then $ m +n \in M$ and $m \in M \cap \frac{1}{p}A = A $, a contradiction.)
Hence \eqref{eq-phi} holds in this case too.
\end{proof}

As a corollary, we recover \cite[Corollary 6.3]{MR2236607}:

\begin{cor}\label{cor-normal-F-pure}
$ k[A]$ is normal, which is of course independent of the characteristic, if and only if it is weakly normal (equivalently F-pure) in arbitrary positive characteristic.
\end{cor}

\begin{proof}
$k[A]$ is normal if and only if for any $m \in M$ and $n \in \Zp$ with $nm \in A$, we have
$m \in A$ if and only if for any $m \in M$ and every prime number $p$ with $pm \in A$, we have
$m \in A$. The last condition is equivalent to that 
$k[A]$ is weakly normal in every positive characteristic.
\end{proof}

\begin{ex}[An example where the monotonicity fails]\label{ex}
Suppose $k$ has characteristic $2$
and  $A \subset \Znn$ is the monoid generated by $8,9,10,11$.
\[
A = \{0,8,9,10,11, 16,17,18,\dots\}
\]
Then the associated 1-dimensional toric variety  $X$
 is not weakly normal, nor F-pure.
 We claim that $\FB_{1} (X)$
is smooth but $\FB_{2} (X)$ is not. In particular $\FB_{2}(X)$ does not dominate
$\FB_{1}(X)$.

For each $e$, $\Delta_{e}$ contains only one 1-dimensional cone,
say $\sigma_{e}$. 
Then
\begin{align*}
B_{\sigma_{1}} & = \frac{1}{2}A\setminus \{ 0 , \frac{9}{2} \} =\{ 4 , 5 , \frac{11}{2} , 8, \frac{17}{2}, 9 \dots \},\\
B_{\sigma_{2}} & =\frac{1}{4}A\setminus \{ 0 , \frac{9}{4},\frac{5}{2},\frac{11}{4} \}
 = \{2, 4,\frac{17}{4},\frac{9}{2}, \dots \} . 
\end{align*}
Since $1 = 11/2 -9/2 \in C _{ \sigma_{1}}$,  $\FB_{1}(X) = \Spec k[x]$.
On the other hand, since none of
\[
 0+1 , \frac{9}{4}+1,\frac{5}{2}+1,\frac{11}{4}+1
\]
 belong to $\frac{1}{4} A$, 
$1 \notin C_{\sigma_{2}}$.
Indeed $C_{\sigma_{2}} = \langle 2,3 \rangle$ 
and  $\FB_{2}(X)=\Spec k[x^{2},x^{3}]$. 
\end{ex}

\begin{ex}[An example where the monotonicity holds, but Condition (*) fails.]
Suppose again $k$ has characteristic $2$
and $A = \langle 2,3\rangle$.
Then for every $e >0$, $\FB_{e}X \cong \Spec k[x]$, because it is the only nontrivial blowup
of $X$. In particular, the monotonicity holds.
We have
\begin{align*}
B_{\sigma_{1}} & = \frac{1}{2}A\setminus \{ 0 , \frac{3}{2} \} =\{ 1,2,\frac{5}{2}, \dots \},\\
B_{\sigma_{2}} & =\frac{1}{4}A\setminus \{ 0 , \frac{1}{2},\frac{3}{4},\frac{5}{4} \}
 = \{1,\frac{3}{2},\frac{7}{4}, \dots  \} . 
\end{align*}
Then $B_{\sigma_{1}} \ne B_{\sigma_{2}} \cap \frac{1}{2}A$. 
Indeed $3/2$ only belongs to the right hand side. 
Therefore Condition (*) in Theorem \ref{thm-local-condition}
fails.
\end{ex}

\section{Stability of F-blowup sequences}\label{sec-stability}

\begin{defn}
Let $X_{1},X_{2},\dots$ be a sequence of blowups of some variety $X$.
Then we say that the sequence \emph{stabilizes} if 
$\exists e_{0}$, $ \forall e \ge e_{0}$, the natural
birational map $X_{e+1} \dasharrow X_{e}$ extends to an isomorphism.

We say that the sequence is \emph{bounded} if
there exists a blowup $Y$ of $X$ which dominates all the $X_{i}$, $i \ge 1$.
\end{defn}

The stability obviously implies the boundedness. 
Conversely, from the following lemma,
 the boundedness together with the monotonicity
implies the stability.

\begin{lem}\label{lem-stability}
Let 
\[
 X_{0} \xleftarrow{f_{1}} X_{1} \xleftarrow{f_{2}} X_{2} \xleftarrow{f_{3}} \cdots
\]
be a sequence of proper surjective morphisms of varieties,
and $g_{i} :Y \to X_{i}$, $i \ge 0$, surjective morphisms of varieties
such that for every $i$, $f_{i} \circ g_{i} = g_{i-1}$.
Then for sufficiently large $i$, $f_{i}$ is an isomorphism.
\end{lem}

\begin{proof}
(Though this fact is perhaps well-known, we
 include a proof for the sake of completeness.) 
Let $\Gamma_{i} \subset Y \times_{k} X_{i}$ be the graph of $g_{i}$
and $H_{i} \subset Y \times_{k} Y$ its inverse image by 
\[
\id_{Y} \times g_{i} :Y \times_{k} Y \to Y \times_{k} X_{i}.
\]
Then we have
\[
H_{i} = \bigsqcup _{y \in Y} \{ y \} \times g_{i}^{-1} (g_{i} (y)) .
\]
Clearly $H_{i -1} \supset H_{i}$. 
Since $Y \times_{k} Y$ has the Noetherian underlying topological space,
for sufficiently large $i$, $H_{i-1}=H_{i}$ and so 
$f_{i}$ is injective and finite.

Now we may suppose that the $f_{i}$ are finite and the $X_{i}$
are affine, say $X_{i} =\Spec R_{i}$. If $S$ denotes the integral closure
of $R_{0}$, then $(R_{i})_{i \in \Znn}$ is an ascending chain of  
$R_{0}$-submodules of $S$. Since $S$ is a Noetherian $R_{0}$-module,
the chain stabilizes.
\end{proof}

\begin{ex}\label{examples}
\begin{enumerate}
\item If $X$ is a 1-dimensional variety, then for sufficiently large $e$,
$\FB_{e}(X)$ is the normalization of $X$ \cite[Corollary 3.19]{Yasuda:math0706.2700}.
In particular, the F-blowup sequence stabilizes. 
\item If $G \subset GL_{d}(k)$ is a finite subgroup of order prime to $p$
and $X = \AA^{d}_{k} /G$, then  for sufficiently large $e$,
$\FB_{e}(X)$ is isomorphic to the $G$-Hilbert scheme $\GHilb (\AA^{d}_{k})$
\cite[Theorem 4.4]{Yasuda:math0706.2700}, \cite[Theorem 1.3]{Yasuda:arXiv:0810.1804}. 
Hence the F-blowup sequence stabilizes. 
\item  The F-blowup sequence of a toric variety is bounded \cite[Theorem 3.13]{Yasuda:math0706.2700}.
Hence the F-blowup sequence of a weakly normal toric variety 
stabilizes. For the normal case, this has been already proved in \cite[Theorem 3.12]{Yasuda:math0706.2700}.
\item Let $R$ be a Noetherian complete local domain over $k$ and $X=\Spec R$.
We can define the $e$-th F-blowup of $X$ as the universal flattening of 
$R^{1/p^{e}}$, see \cite[\S 2.3.2]{Yasuda:math0706.2700}.
We say that  $X$ has
 \emph{finite F-representation type} if there appear only finitely many 
indecomposable $R$-modules, say
$M_{1} , \dots, M_{n}$, up to isomorphism
in $R^{1/p^{i}}$, $i \ge 0$, as direct summands, see \cite[Definition 3.1.1]{MR1444312}. 
If $X$ has finite F-representation type, then the F-blowup sequence of $X$ is bounded. Indeed 
if a blowup $Y \to X$ is a flattening of $N:=\bigoplus_{i=1}^{n} M_{i}$,
then $Y$ dominates all the F-blowups.
Moreover if $X$ is F-pure, then for sufficiently large $e$,
$R^{1/p^{e}}$ has exactly $M_{1} , \dots, M_{n}$
as indecomposable direct summands.
Then $\FB_{e} (X)$ is the universal flattening of $N$.
In particular, the F-blowup sequence stabilizes.  
For instance, every simple singularity  has finite F-representation type.
See \cite{MR1033443} for  simple singularities
in positive characteristic.
Moreover as in the following lemma, every simple singularity of dimension $\ge 3$
is F-pure.
\end{enumerate}
\end{ex}

\begin{lem}
Let $R=k[[x_{0}, \dots, x_{n}]]/(f)$ be a simple hypersurface
singularity of dimension $n \ge 3$. Then $R$ is F-pure.
\end{lem}

\begin{proof}
From the classification \cite{MR1033443}, we may suppose that
$f$ is of the form
\[
f(x_{0}, \dots, x_{n})= g(x_{0}, \dots, x_{n-2})  + x_{n-1}x_{n} .
\]
Then the monomial $x_{n-1}^{p-1}x_{n}^{p-1}$
appears in the expansion of $f^{p-1}$.
Hence $f^{p-1} \notin (x_{0},\dots,x_{n})^{[p]}$. From Fedder's criterion
\cite[Proposition 2.1]{MR701505}, $R$ is F-pure.
\end{proof}

\section{Proof of Proposition \ref{prop-not-preserve}}
\label{sec-not-preserve}

We start with an example of Craw-Maclagan-Thomas \cite[Example 5.7]{MR2356842}.
If $\mathrm{char}(k) \ne 5$, then there exists a finite abelian subgroup 
$G \subset GL_{k} (6)$ of order $5^{4}$ such that the associated
$G$-Hilbert scheme $\GHilb(\AA_{k}^{6})$ is nonnormal. 
Let $X := \AA^{6}_{k}/G$, 
which is a normal  toric (hence F-pure) variety.

As  in Example \ref{examples}, the F-blowup sequence of $X$
stabilizes. For sufficiently large $e$, we have
\[
\FB_{\infty} (X) := \FB_{e} (X) \cong \Hilb^{G}(\AA^{6})
\]
In particular $\FB_{\infty}(X)$ is nonnormal.

But $\FB_{\infty} (X)$ is independent of the base field \cite[Theorem 3.12]{Yasuda:math0706.2700}:
No matter what the base field is, 
the combinatorial data defining the toric variety $\FB_{\infty} (X)$
does not change.
 In particular $\FB_{\infty} (X)$ is well-defined and nonnormal 
also in characteristic $5$. 
From Corollary \ref{cor-normal-F-pure}, it
is not weakly normal nor F-pure in every positive characteristic. 

\bibliographystyle{amsplain}
\bibliography{mybib}

\end{document}